\documentclass[12pt,a4paper]{article} 
 
\usepackage{amsmath,amssymb}
 \usepackage{fullpage}
\usepackage{url}

\newcommand{\dx}{\mathrm{d}} 
\newcommand{\eps}{\varepsilon}
\newcommand{\Q}{\mathbb{Q}}
\newcommand{\R}{\mathbb{R}}
\newcommand{\D}{\mathfrak{D}}
\newcommand{\M}{\mathfrak{M}}
\newcommand{\m}{\mathfrak{m}}
\newcommand{\tr}{\mathfrak{t}}
\newcommand{\singseries}{\mathfrak{S}}
\newcommand{\X}{\mathfrak{X}}
\newcommand{\Odip}[2]{\mathcal{O}_{#1}\!\left(#2\right)\mathchoice{\!}{}{}{}}
\newcommand{\Odi}[1]{\Odip{}{#1}}
\newcommand{\odip}[2]{{o}_{#1}\!\left(#2\right)\mathchoice{\!}{}{}{}}
\newcommand{\odi}[1]{\odip{}{#1}}

\newenvironment{Proof}[1][Proof]{\par\noindent\textbf{#1.}~}
{\hfill$\square$\smallskip\par}
\newtheorem{Theorem}{Theorem}
\newtheorem{Lemma}{Lemma}

\allowdisplaybreaks

\title{A Diophantine problem with a prime and three squares of primes}

\author{Alessandro Languasco and Alessandro Zaccagnini}

\date{}

\begin{document}

\maketitle

\begin{abstract}
We prove that if $\lambda_1$, $\lambda_2$, $\lambda_3$ and $\lambda_4$
are non-zero real numbers, not all of the same sign,
$\lambda_1 / \lambda_2$ is irrational, and $\varpi$ is any real number
then, for any $\eps > 0$ the inequality
$
  \bigl\vert
    \lambda_1 p_1 + \lambda_2 p_2^2 + \lambda_3 p_3^2 + \lambda_4 p_4^2
    +
    \varpi
  \bigr\vert
  \le
  \bigl( \max_j p_j \bigr)^{-1 / 18 + \eps}
$
has infinitely many solution in prime variables $p_1$, \dots, $p_4$.

\noindent
2010 \emph{Mathematics Subject Clas\-sification}: Primary 11D75;
Secon\-da\-ry 11J25, 11P32, 11P55.


\end{abstract}

\section{Introduction}
\label{sec:intro}

This paper deals with an improvement of the recent result of Li and
Wang \cite{LiW2011} concerning Diophantine approximation by means of a
prime and three squares of primes.
We prove the following Theorem.

\begin{Theorem}
\label{Th:main}
Assume that $\lambda_1$, $\lambda_2$, $\lambda_3$ and $\lambda_4$ are
non-zero real numbers, not all of the same sign and that
$\lambda_1 / \lambda_2$ is irrational.
Let $\varpi$ be any real number.
For any $\eps > 0$ the inequality
\begin{equation}
\label{main-ineq}
  \bigl\vert
    \lambda_1 p_1 + \lambda_2 p_2^2 + \lambda_3 p_3^2 + \lambda_4 p_4^2
    +
    \varpi
  \bigr\vert
  \le
  \bigl( \max_j p_j \bigr)^{-1 / 18 + \eps}
\end{equation}
has infinitely many solution in prime variables $p_1$, \dots, $p_4$.
\end{Theorem}

Li and Wang \cite{LiW2011} had $1 / 28$ in place of $1 / 18$.
Our improvement of their result derives from a more efficient use of
Ghosh's bound for exponential sums over squares of primes in
\cite{Ghosh1981} to bound the contribution of the so-called
``intermediate arc.''
This enables us to use a wider ``major arc'' and yields a stronger
result.
The exponent $1 / 18$ arises from there.
We also avoid estimating exponential integrals too early, and we
evaluate them as far as possible, in order to prevent crucial losses
of precision.
We point out that
we can not follow the argument leading to the upper bound for the error term
in formula (3) of \cite{LiW2011}: it does not seem to follow from a
suitable form of the explicit formula by a simple partial integration.
See also the proof of Lemma~5 of Vaughan \cite{Vaughan1974a} or
Lemma~7 of \cite{Vaughan1974b}.

We may change the hypothesis in Theorem~\ref{Th:main} to the
assumption that $\lambda_2 / \lambda_3$ is irrational, say, and the
result is the same, with minor changes in detail.
Furthermore, since the role of $\lambda_2$, $\lambda_3$ and
$\lambda_4$ in our statement above is symmetrical, the assumption that
$\lambda_1 / \lambda_2$ is irrational is not restrictive.

The same kind of argument for the intermediate arc can be used to
improve the result in Languasco and Zaccagnini
\cite{LanguascoZaccagnini2010c}.
For brevity, we simply state the final result, with a very short
sketch of the proof, at the end of this paper.

\section{Outline of the proof}
\label{sec:outline}

We use the variant of the circle method introduced by Davenport and
Heilbronn to deal with Diophantine problems.
In order to prove that \eqref{main-ineq} has infinitely many
solutions, it is sufficient to construct an increasing sequence $X_n$
with limit $+\infty$ such that \eqref{main-ineq} has at least a
solution with $\max_j p_j \in [\delta X_n, X_n]$, where $\delta$ is a
small, fixed positive constant that depends on the coefficients
$\lambda_j$.
This sequence actually depends on rational approximations for
$\lambda_1 / \lambda_2$: more precisely, there are infinitely many
pairs of integers $a$ and $q$ such that $(a, q) = 1$, $q > 0$ and
\[
  \Bigl\vert \frac{\lambda_1}{\lambda_2} - \frac aq \Bigr\vert
  \le
  \frac 1{q^2}.
\]
We take the sequence $X = q^{9 / 5}$ (dropping the useless suffix $n$)
and then, as customary, define all of the circle-method parameters in
terms of $X$.
We may obviously assume that $q$ is sufficiently large.
The choice of the exponent $9 / 5$ is justified in the discussion
following the proof of Lemma~\ref{Lemma-approx}.
Let
\[
  S_1(\alpha)
  =
  \sum_{\delta X \le p \le X} \log p\ e(p \alpha)
  \qquad\text{and}\qquad
  S_2(\alpha)
  =
  \sum_{\delta X \le p^2 \le X} \log p\ e(p^2 \alpha),
\]
where $e(\alpha) = e^{2 \pi i \alpha}$.
As usual, we approximate to $S_1$ and $S_2$ using the functions
\[
  T_1(\alpha)
  =
  \int_{\delta X}^X e(t \alpha) \, \dx t
  \qquad\text{and}\qquad
  T_2(\alpha)
  =
  \int_{(\delta X)^{1/2}}^{X^{1/2}} e(t^2 \alpha) \, \dx t
\]
and notice the simple inequalities
\begin{equation}
\label{bd-Ti}
  T_1(\alpha)
  \ll_\delta
  \min\bigl( X, \vert \alpha \vert^{-1} \bigr)
  \quad\text{and}\quad
  T_2(\alpha)
  \ll_\delta
  X^{-1 / 2}
  \min\bigl( X, \vert \alpha \vert^{-1} \bigr).
\end{equation}
We detect solutions of \eqref{main-ineq} by means of the function
\[
  \widehat{K}_{\eta}(\alpha)
  =
  \max(0, \eta - \vert \alpha \vert)
\]
for $\eta > 0$, which, as the notation suggests, is the Fourier
transform of
\[
  K_{\eta}(\alpha)
  =
  \Bigl( \frac{\sin(\pi \eta \alpha)}{\pi \alpha} \Bigr)^2
\]
for $\alpha \ne 0$, and, by continuity, $K_{\eta}(0) = \eta^2$.
This relation transforms the problem of counting solutions of the
inequality \eqref{main-ineq} into estimating suitable integrals.
We recall the trivial property
\begin{equation}
\label{bd-K(eta)}
  K_{\eta}(\alpha)
  \ll
  \min \Bigl( \eta^2, \vert \alpha \vert^{-2} \Bigr).
\end{equation}

For any measurable subset $\X$ of $\R$ let
\[
  I(\eta, \varpi, \X)
  =
  \int_{\X} S_1(\lambda_1 \alpha) S_2(\lambda_2 \alpha)
    S_2(\lambda_3 \alpha) S_2(\lambda_4 \alpha) K_{\eta}(\alpha)
    e(\varpi \alpha) \, \dx \alpha.
\]
In practice, we take as $\X$ either an interval or a half line, or the
union of two such sets.
The starting point of the method is the observation that
\begin{align*}
  I(\eta, \varpi, \R)
  &=
  \sum_{\substack{\delta X \le p_1 \le X \\ \delta X \le p_j^2 \le X}}
    \log p_1 \log p_2  \log p_3  \log p_4  \\
  &\qquad\qquad\times
  \int_{\R}
    K_{\eta}(\alpha)
    e \bigl(
        (\lambda_1 p_1 + \lambda_2 p_2^2 + \lambda_3 p_3^2 + \lambda_4 p_4^2
        + \varpi) \alpha
      \bigr) \, \dx \alpha \\
  &=
  \sum_{\substack{\delta X \le p_1 \le X \\ \delta X \le p_j^2 \le X}}
    \log p_1 \log p_2  \log p_3 \log p_4   \\
  &\qquad\qquad\times
  \max(0,
       \eta - \vert \lambda_1 p_1 + \lambda_2 p_2^2 + \lambda_3 p_3^2
                + \lambda_4 p_4^2 + \varpi
              \vert) \\
  &\le
  \eta (\log X)^4 \mathcal{N}(X),
\end{align*}
where $\mathcal{N}(X)$ denotes the number of solutions of the
inequality \eqref{main-ineq} with $p_1 \in [\delta X, X]$ and
$p_j^2 \in [\delta X, X]$ for $j = 2$, $3$ and $4$.
We now give the definitions that we need to set up the method.
More definitions will be given at appropriate places later.
We let $P = P(X) = X^{2 / 5} / \log X$,
$\eta = \eta(X) = X^{- 1 / 18 + \eps} (\log X)^2$, and
$R = R(X) = \eta^{-2} (\log X)^2$.
The choice for $P$ is justified at the end of \S\ref{subs:J5}, the one
for $\eta$ at the end of \S\ref{sec:intermediate} and the one for $R$
at the end of \S\ref{sec:trivial}.
We now decompose $\R$ as $\M \cup \m \cup \tr$ where
\[
  \M
  =
  \Bigl[ -\frac PX, \frac PX \Bigr],
  \qquad
  \m
  =
  \Bigl( -R, -\frac PX \Bigr) \cup \Bigl(\frac PX, R \Bigr),
  \qquad
  \tr
  =
  \R \setminus( \M \cup \m),
\]
so that
\[
  I(\eta, \varpi, \R)
  =
  I(\eta, \varpi, \M)
  +
  I(\eta, \varpi, \m)
  +
  I(\eta, \varpi, \tr).
\]
These sets are called the major arc, the intermediate (or minor) arc
and the trivial arc respectively.
In \S\ref{sec:major} we prove that the major arc yields the main term
for $I(\eta, \varpi, \R)$.
In order to show that the contribution of the intermediate arc does
not cancel the main term, we exploit the hypothesis that
$\lambda_1 / \lambda_2$ is irrational to prove that
$\vert S_1(\lambda_1 \alpha) \vert$ and
$\vert S_2(\lambda_2 \alpha) \vert^2$ can not both be large for
$\alpha \in \m$: see \S\ref{sec:intermediate}, and in particular
Lemma~\ref{Lemma-approx}, for the details.
The trivial arc, treated in \S\ref{sec:trivial}, only gives a rather
small contribution.

In the following sections, implicit constants may depend on the
coefficients $\lambda_j$, on $\delta$ and on $\varpi$.

\section{The major arc}
\label{sec:major}

We write
\begin{align*}
  I&(\eta, \varpi, \M)
  =
  \int_{\M}
    S_1(\lambda_1 \alpha) S_2(\lambda_2 \alpha) S_2(\lambda_3 \alpha)
    S_2(\lambda_4 \alpha) K_{\eta}(\alpha) e(\varpi \alpha) \, \dx \alpha \\
  &=
  \int_{\M}
    T_1(\lambda_1 \alpha) T_2(\lambda_2 \alpha) T_2(\lambda_3 \alpha)
    T_2(\lambda_4 \alpha) K_{\eta}(\alpha) e(\varpi \alpha) \, \dx \alpha \\
  &\qquad+
  \int_{\M}
    \bigl( S_1(\lambda_1 \alpha) - T_1(\lambda_1 \alpha) \bigr)
    T_2(\lambda_2 \alpha) T_2(\lambda_3 \alpha) T_2(\lambda_4 \alpha)
    K_{\eta}(\alpha) e(\varpi \alpha) \, \dx \alpha \\
  &\qquad+
  \int_{\M}
    S_1(\lambda_1 \alpha)
    \bigl(S_2(\lambda_2 \alpha) - T_2(\lambda_2 \alpha) \bigr)
    T_2(\lambda_3 \alpha) T_2(\lambda_4 \alpha)
    K_{\eta}(\alpha) e(\varpi \alpha) \, \dx \alpha \\
  &\qquad+
  \int_{\M}
    S_1(\lambda_1 \alpha) S_2(\lambda_2 \alpha)
    \bigl(S_2(\lambda_3 \alpha) - T_2(\lambda_3 \alpha) \bigr)
    T_2(\lambda_4 \alpha)
    K_{\eta}(\alpha) e(\varpi \alpha) \, \dx \alpha \\
  &\qquad+
  \int_{\M}
    S_1(\lambda_1 \alpha) S_2(\lambda_2 \alpha) S_2(\lambda_3 \alpha)
    \bigl(S_2(\lambda_4 \alpha) - T_2(\lambda_4 \alpha) \bigr)
    K_{\eta}(\alpha) e(\varpi \alpha) \, \dx \alpha \\
  &=
  J_1 + J_2 + J_3 + J_4 + J_5,
\end{align*}
say.
We will give a lower bound for $J_1$ and upper bounds for $J_2$,
\dots, $J_5$.
For brevity, since the computations for $J_3$ and $J_4$ are similar
to, but simpler than, the corresponding ones for $J_2$ and $J_5$, we
will skip them.

\subsection{Lower bound for $J_1$}

Apart from very small changes, the lower bound
$J_1 \gg \eta^2 X^{3 / 2}$ is contained in Lemma~8 of Li and Wang
\cite{LiW2011}.
Here we give the required result only in one case, the other ones
being similar.
We have
\begin{align*}
  J_1
  &=
  \int_{\M}
    T_1(\lambda_1 \alpha) T_2(\lambda_2 \alpha) T_2(\lambda_3 \alpha)
    T_2(\lambda_4 \alpha) K_{\eta}(\alpha) e(\varpi \alpha) \, \dx \alpha \\
  &=
  \int_{\R}
    T_1(\lambda_1 \alpha) T_2(\lambda_2 \alpha) T_2(\lambda_3 \alpha)
    T_2(\lambda_4 \alpha) K_{\eta}(\alpha) e(\varpi \alpha) \, \dx \alpha \\
  &\qquad+
  \Odi{
  \int_{P / X}^{+\infty}
    \vert
      T_1(\lambda_1 \alpha) T_2(\lambda_2 \alpha) T_2(\lambda_3 \alpha)
      T_2(\lambda_4 \alpha)
    \vert K_{\eta}(\alpha) \, \dx \alpha}.
\end{align*}
Using inequalities \eqref{bd-Ti} and \eqref{bd-K(eta)}, we see that
the error term is
\[
  \ll
  \eta^2
  X^{-3 / 2}
  \int_{P / X}^{+\infty}
    \frac{\dx \alpha}{\alpha^4}
  \ll
  \eta^2 X^{3 / 2} P^{- 3}
  =
  \odi{\eta^2 X^{3 / 2}}.
\]
For brevity, we set
$\D = [\delta X, X] \times [(\delta X)^{1/2}, X^{1/2}]^3$.
We can rewrite the main term in the form
\begin{align*}
  &
  \idotsint_{\D}
    \int_{\R}
      e \bigl(
      (\lambda_1 t_1 + \lambda_2 t_2^2 + \lambda_3 t_3^2
       +
       \lambda_4 t_4^2 + \varpi) \alpha
        \bigr) \, K_{\eta}(\alpha) \, \dx \alpha
        \, \dx t_1 \, \dx t_2 \, \dx t_3 \, \dx t_4 \\
  &=
  \idotsint_{\D}
    \max(0, \eta
             -
             \vert \lambda_1 t_1 + \lambda_2 t_2^2 + \lambda_3 t_3^2
                   +
                   \lambda_4 t_4^2 + \varpi
             \vert)
    \, \dx t_1 \, \dx t_2 \, \dx t_3 \, \dx t_4.
\end{align*}
We now proceed to show that the last integral is $\gg \eta^2 X^{3 / 2}$.
Apart from trivial changes of sign, there are essentially three cases:

\begin{enumerate}

\item
$\lambda_1 > 0$, $\lambda_2 < 0$,  $\lambda_3 < 0$,  $\lambda_4 < 0$.

\item
$\lambda_1 > 0$, $\lambda_2 > 0$,  $\lambda_3 < 0$,  $\lambda_4 < 0$.

\item
$\lambda_1 > 0$, $\lambda_2 > 0$,  $\lambda_3 > 0$,  $\lambda_4 < 0$.

\end{enumerate}

We briefly deal with the second case.
A suitable change of variables shows that
\begin{align*}
  J_1
  &\gg
  \idotsint_{\D'}
    \max(0, \eta
            -
            \vert \lambda_1 u_1 + \lambda_2 u_2 + \lambda_3 u_3 + \lambda_4 u_4
            \vert)
    \, \frac{\dx u_1 \, \dx u_2 \, \dx u_3 \, \dx u_4}{(u_2 u_3 u_4)^{1/2}} \\
  &\!\!\!\gg
  X^{- 3/2}
  \iiiint_{\D'}\!\!
    \max(0, \eta
            -
            \vert \lambda_1 u_1 + \lambda_2 u_2 + \lambda_3 u_3 + \lambda_4 u_4
            \vert)
    \, \dx u_1 \dx u_2 \dx u_3 \dx u_4,
\end{align*}
where $\D' = [\delta X, (1 - \delta) X]^4$, for large $X$.
For $j = 1$, $2$ and $3$ let
$a_j = 4 \vert \lambda_4 \vert \delta / \vert \lambda_j \vert$,
$b_j = 3 a_j / 2$ and $\mathfrak{I}_j = [a_j X, b_j X]$.
Notice that if $u_j \in \mathfrak{I}_j$ for $j = 1$, $2$ and $3$ then
\[
  \lambda_1 u_1 + \lambda_2 u_2 + \lambda_3 u_3
  \in
  \bigl[
    2 \vert \lambda_4 \vert \delta X, 8 \vert \lambda_4 \vert \delta X
  \bigr]
\]
so that, for every such choice of $(u_1, u_2, u_3)$, the interval
$[a, b]$ with endpoints
$\pm \eta / \vert \lambda_4 \vert +
(\lambda_1 u_1 + \lambda_2 u_2 + \lambda_3 u_3) / \vert \lambda_4 \vert$
is contained in $[\delta X, (1 - \delta) X]$.
In other words, for $u_4 \in [a, b]$ the values of
$\lambda_1 u_1 + \lambda_2 u_2 + \lambda_3 u_3 + \lambda_4 u_4$
cover the whole interval $[-\eta, \eta]$.
Hence, for any $(u_1, u_2, u_3) \in
\mathfrak{I}_1 \times \mathfrak{I}_2 \times \mathfrak{I}_3$ we have
\begin{align*}
  \int_{\delta X}^{(1 - \delta) X}
   & \max(0, \eta
            -
            \vert \lambda_1 u_1 + \lambda_2 u_2 + \lambda_3 u_3 + \lambda_4 u_4
            \vert)
    \, \dx u_4 
    \\
    &
  =
  \vert \lambda_4 \vert^{-1}
  \int_{-\eta}^{\eta} \max(0, \eta - \vert u \vert) \, \dx u
  \gg
  \eta^2.
\end{align*}
Finally,
\[
  J_1
  \gg
  \eta^2
  X^{- 3 / 2}
  \iiint_{\mathfrak{I}_1 \times \mathfrak{I}_2 \times \mathfrak{I}_3}
    \dx u_1 \, \dx u_2 \, \dx u_3
  \gg
  \eta^2 X^{3 / 2},
\]
which is the required lower bound.

\subsection{Bound for $J_2$}

Let
\[
  U_1(\alpha)
  =
  \sum_{\delta X \le n \le X} e(n \alpha)
  \qquad\text{and}\qquad
  U_2(\alpha)
  =
  \sum_{\delta X \le n^2 \le X} e(n^2 \alpha).
\]
By the Euler summation formula we have
\begin{equation}
\label{bd-T-U}
  T_j(\alpha)
  -
  U_j(\alpha)
  \ll
  1 + \vert \alpha \vert X
  \qquad\text{for $j = 1$, $2$.}
\end{equation}
Using \eqref{bd-K(eta)} we see that
\begin{align*}
  J_2
  &\ll
  \eta^2
  \int_{\M}
    \bigl\vert S_1(\lambda_1 \alpha) - T_1(\lambda_1 \alpha) \bigr\vert \,
    \vert T_2(\lambda_2 \alpha) \vert \,
    \vert T_2(\lambda_3 \alpha) \vert \,
    \vert T_2(\lambda_4 \alpha) \vert \, \dx \alpha \\
  &\le
  \eta^2
  \int_{\M}
    \bigl\vert S_1(\lambda_1 \alpha) - U_1(\lambda_1 \alpha) \bigr\vert \,
    \vert T_2(\lambda_2 \alpha) \vert \,
    \vert T_2(\lambda_3 \alpha) \vert \,
    \vert T_2(\lambda_4 \alpha) \vert \, \dx \alpha \\
  &\qquad+
  \eta^2
  \int_{\M}
    \bigl\vert U_1(\lambda_1 \alpha) - T_1(\lambda_1 \alpha) \bigr\vert \,
    \vert T_2(\lambda_2 \alpha) \vert \,
    \vert T_2(\lambda_3 \alpha) \vert \,
    \vert T_2(\lambda_4 \alpha) \vert \, \dx \alpha \\
  &=
  \eta^2
  (A_2 + B_2),
\end{align*}
say.
In order to estimate $A_2$ we connect it to the Selberg integral as in
Lemma~6 of Languasco and Zaccagnini \cite{LanguascoZaccagnini2010c}.
We set
\[
  J(X, h)
  =
  \int_{\delta X}^X \bigl( \theta(x + h) - \theta(x) - h)^2 \, \dx x,
\]
where $\theta$ is the usual Chebyshev function.
By the Cauchy inequality and \eqref{bd-Ti} above, for any fixed
$A > 0$ we have
\begin{align*}
  A_2
  &\ll
  \Bigl(
    \int_{-P / X}^{P / X}
      \bigl\vert S_1(\lambda_1 \alpha) - U_1(\lambda_1 \alpha) \bigr\vert^2
      \, \dx \alpha
  \Bigr)^{1 / 2}
  \\
  &
  \hskip 1cm\times
  \Bigl(
    \int_{-P / X}^{P / X}
      \vert T_2(\lambda_2 \alpha) \vert^2 \,
      \vert T_2(\lambda_3 \alpha) \vert^2 \,
      \vert T_2(\lambda_4 \alpha) \vert^2 \, \dx \alpha
  \Bigr)^{1 / 2} \\
  &\ll
  \frac PX J \Bigl(X, \frac XP \Bigr)^{1 / 2}
  \Bigl(
    \int_0^{1 / X} X^3 \, \dx \alpha
    +
    \int_{1 / X}^{P / X} \frac{\dx \alpha}{X^3 \alpha^6}
  \Bigr)^{1 / 2} \\
  &\ll_A
  \Bigl(\frac X{(\log X)^A} \Bigr)^{1 / 2}
  X
  \ll_A
  \frac{X^{3 / 2}}{(\log X)^{A / 2}}
\end{align*}
by the Theorem in \S6 of Saffari and Vaughan \cite{SaffariV1977},
which we can use provided that $X / P \ge X^{1 / 6 + \eps}$, that is,
$P \le X^{5 / 6 - \eps}$.
This proves that $\eta^2 A_2 = \odi{\eta^2 X^{3 / 2}}$.
Furthermore, using the inequalities \eqref{bd-Ti} and \eqref{bd-T-U}
we see that
\begin{align*}
  B_2
  &\ll
  \int_0^{1 / X}
    \vert T_2(\lambda_2 \alpha) \vert \,
    \vert T_2(\lambda_3 \alpha) \vert \,
    \vert T_2(\lambda_4 \alpha) \vert \, \dx \alpha
    \\
    &
    \hskip1cm
  +
  X
  \int_{1 / X}^{P / X}
    \alpha \,
    \vert T_2(\lambda_2 \alpha) \vert \,
    \vert T_2(\lambda_3 \alpha) \vert \,
    \vert T_2(\lambda_4 \alpha) \vert \, \dx \alpha \\
  &\ll
  \frac1X X^{3 / 2}
  +
  X
  \int_{1 / X}^{P / X} \alpha X^{- 3 / 2} \, \frac{\dx \alpha}{\alpha^3}
  \ll
  X^{1 / 2}
  +
  X^{- 1/ 2}
  \int_{1 / X}^{P / X} \frac{\dx \alpha}{\alpha^2}
  \ll
  X^{1 / 2},
\end{align*}
so that $\eta^2 B_2 = \odi{\eta^2 X^{3 / 2}}$.

\subsection{Bound for $J_5$}
\label{subs:J5}

Inequality \eqref{bd-K(eta)} implies that
\begin{align*}
  J_5
  &\ll
  \eta^2
  \int_{\M}
    \bigl\vert S_1(\lambda_1 \alpha) \bigr\vert \,
    \bigl\vert S_2(\lambda_2 \alpha) \bigr\vert \,
    \bigl\vert S_2(\lambda_3 \alpha) \bigr\vert \,
    \bigl\vert S_2(\lambda_4 \alpha) - T_2(\lambda_4 \alpha) \bigr\vert
    \, \dx \alpha \\
  &\ll
  \eta^2
  \int_{\M}
    \bigl\vert S_1(\lambda_1 \alpha) \bigr\vert \,
    \bigl\vert S_2(\lambda_2 \alpha) \bigr\vert \,
    \bigl\vert S_2(\lambda_3 \alpha) \bigr\vert \,
    \bigl\vert S_2(\lambda_4 \alpha) - U_2(\lambda_4 \alpha) \bigr\vert
    \, \dx \alpha \\
  &\quad+
  \eta^2
  \int_{\M}
    \bigl\vert S_1(\lambda_1 \alpha) \bigr\vert \,
    \bigl\vert S_2(\lambda_2 \alpha) \bigr\vert \,
    \bigl\vert S_2(\lambda_3 \alpha) \bigr\vert \,
    \bigl\vert U_2(\lambda_4 \alpha) - T_2(\lambda_4 \alpha) \bigr\vert
    \, \dx \alpha \\
  &=
  \eta^2 (A_5 + B_5),
\end{align*}
say.
Now let
\[
  J^*(X, h)
  =
  \int_{\delta X}^X
    \bigl(
      \theta(\sqrt{x + h}) - \theta(\sqrt{x})
      -
      (\sqrt{x + h} - \sqrt{x})
    \bigr)^2 \, \dx x.
\]
The Parseval inequality and trivial bounds yield, for any fixed
$A > 0$,
\begin{align*}
  A_5
  &\ll
  X
  \Bigl(
    \int_{\M} \bigl\vert S_1(\lambda_1 \alpha) \bigr\vert^2 \, \dx \alpha
  \Bigr)^{1 / 2}
  \Bigl(
    \int_{\M}
      \bigl\vert S_2(\lambda_4 \alpha) - U_2(\lambda_4 \alpha) \bigr\vert^2
      \, \dx \alpha
  \Bigr)^{1 / 2} \\
  &\ll
  X (X \log X)^{1 / 2}
  \frac PX J^*\Bigl(X, \frac XP \Bigr)^{1 / 2} \\
  &\ll_A
  X^{3 / 2} (\log X)^{1 / 2 - A / 2}
\end{align*}
by Lemmas 3.12 and 3.13 of Languasco and Settimi \cite{LanguascoS2012},
which we can use provided that $X / P \ge X^{7 / 12 + \eps}$, that is,
$P \le X^{5 / 12 - \eps}$.
This proves that $\eta^2 A_5 = \odi{\eta^2 X^{3 / 2}}$.
Furthermore, using \eqref{bd-T-U}, the Cauchy inequality and trivial
bounds we see that
\begin{align*}
  B_5
  &\ll
  \int_0^{1 / X}
    \bigl\vert S_1(\lambda_1 \alpha) \bigr\vert \,
    \bigl\vert S_2(\lambda_2 \alpha) \bigr\vert \,
    \bigl\vert S_2(\lambda_3 \alpha) \bigr\vert \, \dx \alpha \\
  &\qquad+
  X
  \int_{1 / X}^{P / X}
    \alpha
    \bigl\vert S_1(\lambda_1 \alpha) \bigr\vert \,
    \bigl\vert S_2(\lambda_2 \alpha) \bigr\vert \,
    \bigl\vert S_2(\lambda_3 \alpha) \bigr\vert \, \dx \alpha \\
  &\ll
  \frac1X X^2
  +
  X
  \Bigl(
    \int_{1 / X}^{P / X}
      \alpha^4 \, \dx \alpha
  \Bigr)^{1 / 4}
  \Bigl(
    \int_{1 / X}^{P / X}
      \bigl\vert S_1(\lambda_1 \alpha) \bigr\vert^2 \, \dx \alpha
  \Bigr)^{1 / 2}  \\
  &\qquad\times
  \max_{\alpha \in [1 / X, P / X]}
    \bigl\vert S_2(\lambda_2 \alpha) \bigr\vert
  \Bigl(
    \int_{1 / X}^{P / X}
    \bigl\vert S_2(\lambda_3 \alpha) \bigr\vert^4 \, \dx \alpha
  \Bigr)^{1 / 4} \\
  &\ll
  X
  +
  X
  \Bigl( \frac PX \Bigr)^{5 / 4}
  (X \log X)^{1 / 2}
  \max_{\alpha \in [1 / X, P / X]}
    \bigl\vert S_2(\lambda_2 \alpha) \bigr\vert
    \\
    &
    \qquad
    \times
  \Bigl(
    \int_0^1
    \bigl\vert S_2(\lambda_3 \alpha) \bigr\vert^4 \, \dx \alpha
  \Bigr)^{1 / 4} \\
  &\ll
  X
  +
  X^{3 / 4} P^{5 / 4} (\log X)^{1 / 2}
  \Bigl(
    \int_0^1
    \bigl\vert S_2(\lambda_3 \alpha) \bigr\vert^4 \, \dx \alpha
  \Bigr)^{1 / 4}.
\end{align*}
In order to estimate the integral at the far right we borrow (4.7)
from Languasco and Settimi \cite{LanguascoS2012}, that gives the bound
$\ll X (\log X)^2$.

\noindent
Hence $B_5$ $\ll $ $X P^{5 / 4} \log X$, so that
$\eta^2 B_5 = \odi{\eta^2 X^{3 / 2}}$ provided that
$P$ $=$ $\odi{X^{2 / 5} (\log X)^{-4 / 5}}$.
We may therefore choose $P = X^{2 / 5} / (\log X)$.

\section{The intermediate arc}
\label{sec:intermediate}

We need to show that $\vert S_1(\lambda_1 \alpha) \vert$ and
$\vert S_2(\lambda_2 \alpha) \vert^2$ can not both be large for
$\alpha \in \m$, exploiting the fact that $\lambda_1 / \lambda_2$ is
irrational.
We do this using two famous results by Vaughan about $S_{1}(\alpha)$
and by Ghosh about $S_{2}(\alpha)$.

\begin{Lemma}[Vaughan \cite{Vaughan1997}, Theorem 3.1]
\label{Vaughan-estim}
Let $\alpha$ be a real number and $a,q$ be positive integers
satisfying $(a,q)=1$ and $\vert \alpha -a/q \vert < q^{-2}$.
Then
\[
  S_{1}(\alpha)
  \ll
  \Bigl(\frac{X}{\sqrt{q}} + \sqrt{Xq} + X^{4/5} \Bigr) \log^4 X.
\]
\end{Lemma}
\begin{Lemma}[Ghosh \cite{Ghosh1981}, Theorem 2]
\label{Ghosh-estim}
Let $\alpha$ be a real number and $a,q$ be positive integers
satisfying $(a,q)=1$ and $\vert \alpha -a/q \vert < q^{-2}$.
Let moreover $\epsilon>0$.
Then
\[
S_{2}(\alpha)
\ll_{\epsilon}
X^{1/2+\epsilon}
\left(
\frac{1}{q}
+
\frac{1}{X^{1/4}}
+
\frac{q}{X}
\right)^{1/4}.
\]
\end{Lemma}

\begin{Lemma}
\label{Lemma-approx}
Assume that $\lambda_1 / \lambda_2$ is irrational and let $X = q^{9 / 5}$,
where $q$ is the denominator of a convergent of the continued fraction
for $\lambda_1 / \lambda_2$.
Let $V(\alpha) =
\min \bigl( \vert S_1(\lambda_1 \alpha) \vert^{1 / 2}$,
           $\vert S_2(\lambda_2 \alpha) \vert \bigr)$.
Then, for arbitrary $\eps > 0$, we have
\[
  \sup_{\alpha \in \m} V(\alpha)
  \ll
  X^{4 / 9 + \eps}.
\]
\end{Lemma}

\begin{Proof}
Let $\alpha \in \m$ and $Q = X^{2/9} / \log X \leq P$.
By Dirichlet's Theorem, there exist integers $a_{i},q_{i}$  with
$1\leq q_{i}\leq X/Q$ and $(a_{i},q_{i})=1$, such that
$\vert \lambda_{i} \alpha q_{i}-a_{i}\vert \leq Q/X$, for $i=1,2$.
We remark that $a_{1}a_{2} \neq 0$ otherwise we would have $\alpha\in \M$.
Now suppose that  $q_{i} \leq Q$ for $i=1,2$. In this case we get
\[
a_{2}q_{1} \frac{\lambda_{1}}{\lambda_{2}} - a_{1}q_{2}
=
( \lambda_{1} \alpha q_{1}-a_{1}) \frac{a_{2}}{\lambda_{2} \alpha}
-
( \lambda_{2} \alpha q_{2}-a_{2}) \frac{a_{1}}{\lambda_{2} \alpha}
\]
and hence
\begin{equation}
\label{bd-1}
\left\vert
a_{2}q_{1} \frac{\lambda_{1}}{\lambda_{2}} - a_{1}q_{2}
\right\vert
\leq
2\left(
1+ \left\vert  \frac{\lambda_{1}}{\lambda_{2}} \right\vert
\right)
\frac{Q^{2}}{X}
<
\frac{1}{2q}
\end{equation}
for sufficiently large $X$.
Then, from the law of best approximation and the definition of $\m$, we obtain
\begin{equation}
\label{bd-2}
  X^{5/9}=q
  \leq
  \vert a_{2}q_{1} \vert
  \ll q_{1}q_{2} R
  \leq Q^{2} R
  \leq X^{5 / 9 - 2 \eps} \log^{-4} X,
\end{equation}
which is absurd.
Hence either $q_{1}>Q$ or $q_{2}>Q$.
Assume first that $q_{2}>Q$.
Using Lemma \ref{Ghosh-estim} on $S_2(\lambda_2 \alpha)$, we have
\begin{align}
\notag
V(\alpha)
\leq
\vert
S_2(\lambda_2 \alpha)
\vert
& \ll_{\eps}
X^{1/2+\eps}
\sup_{Q<q_{2}\leq X/Q }
\left(
\frac{1}{q_{2}}
+
\frac{1}{X^{1/4}}
+
\frac{q_{2}}{X}
\right)^{1/4}
\\
\label{first-minor}
&
\ll_{\eps}
X^{4/9+\eps}
(\log X)^{1/4}.
\end{align}
Assume now that $q_{1}>Q$.
Using Lemma \ref{Vaughan-estim} on $S_1(\lambda_1 \alpha)$, we have
\begin{align}
\notag
V(\alpha)
\leq
\vert
S_1(\lambda_1 \alpha)
\vert^{1/2}
& \ll
\sup_{Q<q_{1}\leq X/Q }
\left(
\frac{X}{\sqrt{q_{1}}}
+
\sqrt{Xq_{1}}
+
X^{4/5}
\right)^{1/2}
\log^{2}X
\\
\label{second-minor}
&
\ll
X^{4/9}
(\log X)^{3}.
\end{align}
Lemma \ref{Lemma-approx} follows combining \eqref{first-minor} and
\eqref{second-minor}.
\end{Proof}

The constraint on the choice $X = q^{9 / 5}$ arises from the bounds
\eqref{bd-1} and \eqref{bd-2}.
Their combination prevents us from choosing the optimal value
$X = q^2$.

\begin{Lemma}
\label{Lemma:bd-minor}
We have
\[
  \int_{\m}
    \vert S_1(\lambda_1 \alpha) \vert^2 K_{\eta}(\alpha) \, \dx \alpha
  \ll
  \eta X \log X
  \]
 and
  \[
  \int_{\m}
    \vert S_2(\lambda_j \alpha) \vert^4 K_{\eta}(\alpha) \, \dx \alpha
  \ll
  \eta X (\log X)^2
\]
for $j = 2$, $3$ and $4$.
\end{Lemma}

\begin{Proof}
The proof is achieved arguing as in \S\ref{sec:trivial} below where we
bound the quantities $A$ and $B$, the main difference being the fact
that we have to split the range $[P / X, R]$ into two intervals in
order to use \eqref{bd-K(eta)} efficiently.
See also the proof of Lemma~12 of \cite{LiW2011}.
For the sake of brevity we skip the details.
\end{Proof}

Now let
\begin{align*}
  \X_1
  &=
  \{ \alpha \in [P / X, R] \colon
    \vert S_1(\lambda_1 \alpha) \vert^{1 / 2}
    \le
    \vert S_2(\lambda_2 \alpha) \vert \} \\
  \X_2
  &=
  \{ \alpha \in [P / X, R] \colon
    \vert S_1(\lambda_1 \alpha) \vert^{1 / 2}
    \ge
    \vert S_2(\lambda_2 \alpha) \vert \}
\end{align*}
so that $[P / X, R] = \X_1 \cup \X_2$ and
\[
  \Bigl\vert I(\eta, \varpi, \m) \Bigr\vert
  \ll
  \Bigl( \int_{\X_1} + \int_{\X_2} \Bigr)
    \bigl\vert
      S_1(\lambda_1 \alpha) S_2(\lambda_2 \alpha)
      S_2(\lambda_3 \alpha) S_2(\lambda_4 \alpha)
    \bigr\vert
    K_{\eta}(\alpha) \, \dx \alpha.
\]
H\"older's inequality gives
\begin{align*}
  \int_{\X_1}
  &\le
  \Bigl(
    \int_{\X_1}
      \vert S_1(\lambda_1 \alpha) \vert^4 K_{\eta}(\alpha) \, \dx \alpha
  \Bigr)^{1 / 4}
  \prod_{j = 2}^4
  \Bigl(
    \int_{\X_1}
      \vert S_2(\lambda_j \alpha) \vert^4 K_{\eta}(\alpha) \, \dx \alpha
  \Bigr)^{1 / 4} \\
  &\le
  \max_{\alpha \in \X_1} \vert S_1(\lambda_1 \alpha) \vert^{1 / 2}
  \Bigl(
    \int_{\m}
      \vert S_1(\lambda_1 \alpha) \vert^2 K_{\eta}(\alpha) \, \dx \alpha
  \Bigr)^{1 / 4}
  \\
  &
  \qquad
  \times
  \prod_{j = 2}^4
  \Bigl(
    \int_{\m}
      \vert S_2(\lambda_j \alpha) \vert^4 K_{\eta}(\alpha) \, \dx \alpha
  \Bigr)^{1 / 4} \\
  &\ll
  X^{4 / 9 + \eps} (\eta X \log X)^{1 / 4} (\eta X (\log X)^2)^{3 / 4} \\
  &\ll
  \eta X^{13 / 9 + \eps} (\log X)^{7 / 4}
\end{align*}
by Lemmas~\ref{Lemma-approx} and \ref{Lemma:bd-minor}.
The computation on $\X_2$ is similar: we have
\begin{align*}
  \int_{\X_2}
  &\le
  \Bigl(
    \int_{\X_2}
      \vert S_1(\lambda_1 \alpha) \vert^2 K_{\eta}(\alpha) \, \dx \alpha
  \Bigr)^{1 / 2}
  \max_{\alpha \in \X_2} \vert S_2(\lambda_2 \alpha) \vert
  \\
  &
  \qquad
  \times
  \prod_{j = 3}^4
  \Bigl(
    \int_{\X_2}
      \vert S_2(\lambda_j \alpha) \vert^4 K_{\eta}(\alpha) \, \dx \alpha
  \Bigr)^{1 / 4} \\
  &\ll
  (\eta X \log X)^{1 / 2} X^{4 / 9 + \eps} (\eta X (\log X)^2)^{1 / 2} \\
  &\ll
  \eta X^{13 / 9 + \eps} (\log X)^{3 / 2},
\end{align*}
again by Lemmas~\ref{Lemma-approx} and \ref{Lemma:bd-minor}.
Summing up,
\[
  \Bigl\vert I(\eta, \varpi, \m) \Bigr\vert
  \ll
  \eta X^{13 / 9 + \eps} (\log X)^{7 / 4},
\]
and this is $\odi{\eta^2 X^{3 / 2}}$ provided that
$\eta \ge X^{- 1 / 18 + \eps} (\log X)^2$.

\section{The trivial arc}
\label{sec:trivial}

Using the Cauchy inequality and a trivial bound for
$S_2(\lambda_4 \alpha)$ we see that
\begin{align*}
  \Bigl\vert &I(\eta, \varpi, \tr) \Bigr\vert
  \le
  2 
  \int_R^{+\infty}
    \vert S_1(\lambda_1 \alpha) \vert \, \vert S_2(\lambda_2 \alpha) \vert \,
    \vert S_2(\lambda_3 \alpha) \vert \, \vert S_2(\lambda_4 \alpha) \vert \,
    K_{\eta}(\alpha) \, \dx \alpha \\
  &\ll
  \sup_{\alpha \in (R, +\infty)} \vert S_2(\lambda_4 \alpha) \vert 
  \Bigl(
    \int_R^{+\infty}
      \vert S_1(\lambda_1 \alpha) \vert^2 K_{\eta}(\alpha) \, \dx \alpha
  \Bigr)^{1/2}
   \\
  &\qquad\times
  \Bigl(
    \int_R^{+\infty}
    \vert S_2(\lambda_2 \alpha) \vert^2 \,
    \vert S_2(\lambda_3 \alpha) \vert^2 \, K_{\eta}(\alpha) \, \dx \alpha
  \Bigr)^{1/2}
\\
  &\ll
  X^{1 / 2}
  \Bigl(
    \int_R^{+\infty}
      \vert S_1(\lambda_1 \alpha) \vert^2 K_{\eta}(\alpha) \, \dx \alpha
  \Bigr)^{1/2}
  \Bigl(
    \int_R^{+\infty}
    \vert S_2(\lambda_2 \alpha) \vert^4 \, K_{\eta}(\alpha) \, \dx \alpha
  \Bigr)^{1/4}
   \\
  &\qquad\times
  \Bigl(
    \int_R^{+\infty}
    \vert S_2(\lambda_3 \alpha) \vert^4 \, K_{\eta}(\alpha) \, \dx \alpha
  \Bigr)^{1/4} \\
  &\ll
  X^{1 / 2}
  \Bigl(
    \int_{\vert \lambda_1 \vert R}^{+\infty}
      \frac{\vert S_1(\alpha) \vert^2}{\alpha^2} \, \dx \alpha
  \Bigr)^{1/2}
  \Bigl(
    \int_{\vert \lambda_2 \vert R}^{+\infty}
      \frac{\vert S_2(\alpha) \vert^4}{\alpha^2} \, \dx \alpha
  \Bigr)^{1/4}
  \\
  &
  \qquad \times
  \Bigl(
    \int_{\vert \lambda_3 \vert R}^{+\infty}
      \frac{\vert S_2(\alpha) \vert^4}{\alpha^2} \, \dx \alpha
  \Bigr)^{1/4} \\
  &\ll
  X^{1 / 2}
  A^{1 / 2}
  B^{1 / 2},
\end{align*}
say, where in the last but one line we used the inequality
\eqref{bd-K(eta)}, and we set
\[
  A
  =
  \int_{\vert \lambda_1 \vert R}^{+\infty}
    \frac{\vert S_1(\alpha) \vert^2}{\alpha^2} \, \dx \alpha
  \qquad\text{and}\qquad
  B
  =
  \int_{\min(\vert \lambda_2 \vert, \vert \lambda_3 \vert) R}^{+\infty}
    \frac{\vert S_2(\alpha) \vert^4}{\alpha^2} \, \dx \alpha.
\]
Using periodicity we have
\[
  A
  \ll
  \sum_{n \ge \vert \lambda_1 \vert R}
    \frac 1{(n - 1)^2}
    \int_{n - 1}^n \vert S_1(\alpha) \vert^2 \, \dx \alpha
  \ll
  \frac{X \log X}{\vert \lambda_1 \vert R}
\]
by the Prime Number Theorem, while
\[
  B
  \ll
  \sum_{n \ge \min(\vert \lambda_2 \vert, \vert \lambda_3 \vert) R}
    \frac 1{(n - 1)^2}
    \int_{n-1}^n \vert S_2(\alpha) \vert^4 \, \dx \alpha
  \ll
  \frac{X (\log X)^2}{\min(\vert \lambda_2 \vert, \vert \lambda_3 \vert) R}.
\]
The last estimate follows from Satz~3 of Rieger \cite{Rieger1968},
which is used to bound ``non-diagonal'' solutions of
$p_1^2 + p_2^2 = p_3^2 + p_4^2$, and the Prime Number Theorem for the
remaining solutions.
See also the bound for $H_{12}$ in Liu \cite{Liu2004}.
Collecting these estimates, we conclude that
\begin{equation}
\label{final-est-trivial}
  \Bigl\vert I(\eta, \varpi, \tr) \Bigr\vert
  \ll
  \frac{X^{3 / 2} (\log X)^{3 / 2}}R.
\end{equation}
Hence, the choice $R = \eta^{-2} (\log X)^2$ is admissible.

\section{Proof of Theorem~\ref{Th-LZ}}

In our paper \cite{LanguascoZaccagnini2010c} we dealt with a similar
problem, with two primes and $s$ powers of $2$.
The goal was to approximate any real number by means of values of the
form
\begin{equation}
\label{two-primes}
  \lambda_1 p_1 + \lambda_2 p_2
  +
  \mu_1 2^{m_1} + \dots + \mu_s 2^{m_s},
\end{equation}
where $\lambda_1$ and $\lambda_2$ are real numbers of opposite sign,
with an irrational ratio, and the non-zero coefficients $\mu_1$,
\dots, $\mu_s$ satisfy suitable conditions, $p_1$ and $p_2$ are prime
numbers and $m_1$, \dots, $m_s$ are positive integers.
The result is an upper bound on the least value $s_0$ that ensures the
existence of an approximation of the form \eqref{two-primes} for all
$s \ge s_0$.
The quality of the result depends on rational approximations to
$\lambda_1 / \lambda_2$: we let $\mathfrak{R}$ denote the set of
irrational numbers $\xi$ such that the denominators $q_m$ of the
convergents to $\xi$, arranged in increasing order of magnitude,
satisfy $q_{m + 1} \ll q_m^{1 + \eps}$.
By Roth's Theorem, all algebraic numbers belong to $\mathfrak{R}$, and
almost all real numbers, in the sense of the Lebesgue measure, also
belong to $\mathfrak{R}$.
We denote by $\mathfrak{R}'$ the set of irrational numbers that do not
belong to $\mathfrak{R}$.
For $\lambda_1 / \lambda_2$ belonging to this set, we have the
following improvement of our result in \cite{LanguascoZaccagnini2010c}.

\begin{Theorem}
\label{Th-LZ}
Suppose that $\lambda_1$ and $\lambda_2$ are real numbers such that
$\lambda_1 / \lambda_2$ is negative and irrational with $\lambda_1 > 1$,
$\lambda_2 < -1$ and $\vert \lambda_1 / \lambda_2 \vert \ge 1$.
Further suppose that $\mu_1$, \dots, $\mu_s$ are nonzero real numbers
such that $\lambda_i / \mu_i \in \Q$ for $i \in \{ 1$, $2\}$, and
denote by $a_i / q_i$ their reduced representations as rational
numbers.
Let moreover $\eta$ be a sufficiently small positive constant such
that $\eta < \min(\lambda_1 / a_1; \vert \lambda_2 / a_2 \vert)$.
Finally, for $\lambda_1 / \lambda_2 \in \mathfrak{R}'$, let
\[
  s_0
  =
  2
  +
  \Bigl\lceil
    \frac{\log (C(q_1, q_2) \lambda_1) - \log \eta}
         {-\log (0.884472132)}
  \Bigr\rceil.
\]
Then for every real number $\gamma$ and every integer $s \ge s_0$ the
inequality
\[
  \vert
  \lambda_1 p_1 + \lambda_2 p_2
  +
  \mu_1 2^{m_1} + \dots + \mu_s 2^{m_s}
  +
  \gamma
  \vert
  <
  \eta
\]
has infinitely many solutions in primes $p_1$,$ p_2$ and positive
integers $m_1$, \dots, $m_s$, where
$C(q_1, q_2) =
\Bigl( \log{2} + C \cdot \singseries^{\prime}(q_1) \Bigr)^{1/2}
\Bigl( \log{2} + C \cdot \singseries^{\prime}(q_2) \Bigr)^{1/2}$,
$C = 10.0219168340$ and
\[
  \singseries^{\prime}(n)
  =
  \prod_{\substack{p \mid n \\ p > 2}} \frac{p - 1}{p - 2}.
\]
\end{Theorem}

We can improve our previous treatment of the intermediate arc in \S7
of \cite{LanguascoZaccagnini2010c}.
We let $V(\alpha) = \min \bigl( \vert S_{1}(\lambda_1 \alpha) \vert,
\vert S_{1}(\lambda_2 \alpha) \vert \bigr)$ and recall that $\m_2$ is the
subset of $[X^{-2 / 3}, (\log X)^2]$ where the exponential sum
$G(\alpha) = \sum_{n \le L} e(2^n \alpha)$ is ``large'' in absolute
value.
Here $L = (\log (\eps X / 2 M)) / \log 2$ where
$M = \max_j \vert \mu_j \vert$.
The technique due to Pintz and Ruzsa \cite{PintzR2003} ensures that
its measure is comparatively small.
In the following computation, implicit constants may depend on
$\lambda_1$ and $\lambda_2$.
We have
\begin{align*}
  \Bigl\vert
    \int_{\m_2}
      &S_1(\lambda_1 \alpha) S_1(\lambda_2 \alpha)
      \prod_{j = 1}^s G(\mu_j \alpha)
      K_{\eta}(\alpha) \, \dx \alpha
        \Bigr\vert
      \\
      &
  \ll
  \eta^2 (\log X)^s
  \int_{\m_2}
    \vert S_1(\lambda_1 \alpha) S_1(\lambda_2 \alpha) \vert \, \dx \alpha \\
  &\ll
  \eta^2 (\log X)^s
  \sup_{\alpha \in \m_2} V(\alpha)
  \int_{\m_2}
    \vert S_1(\lambda_2 \alpha) \vert \, \dx \alpha \\
  &\ll
  \eta^2 (\log X)^s
  \sup_{\alpha \in \m_2} V(\alpha)
  \Bigl( \int_{\m_2} \dx \alpha \Bigr)^{1 / 2}
  \Bigl( \int_{\m_2} \vert S_1(\lambda_2 \alpha) \vert^2 \, \dx \alpha
  \Bigr)^{1 / 2} \\
  &\ll
  \eta^2 (\log X)^s
  \vert \m_2 \vert^{1 / 2}
  \bigl( X (\log X)^3 \bigr)^{1 / 2}
  \sup_{\alpha \in \m_2} V(\alpha) \\
  &\ll
  \eta^2 (\log X)^s
  (\log X) s^{1 / 2} X^{-c / 2}
  X^{1 / 2} (\log X)^{3 / 2}
  \sup_{\alpha \in \m_2} V(\alpha) \\
  &\ll
  \eta^2 s^{1 / 2}
  X^{1 / 2 - c / 2}
  (\log X)^{s + 5 / 2}
  \sup_{\alpha \in \m_2} V(\alpha).
\end{align*}
The proof of Lemma~4 of Parsell \cite{Parsell2003} implies that
\[
  \sup_{\alpha \in \m_2} V(\alpha)
  =
  \sup_{\alpha \in \m_2}
    \min
      \bigl(
        \vert S_1(\lambda_1 \alpha) \vert, \vert S_1(\lambda_2 \alpha) \vert
      \bigr)
  \ll
  X^{7 / 8} (\log X)^5.
\]
Hence the integral above is bounded by
\[
  \eta^2 s^{1 / 2}
  X^{11 / 8 - c / 2}
  (\log X)^{s + 15 / 2}
\]
It is therefore sufficient to take $c > \frac 34$ (instead of the
bound $c > \frac 45$ that we had in \cite{LanguascoZaccagnini2010c}).
Taking $c = \frac34 + 10^{-20}$, the method due to Pintz and Ruzsa
(see for example Lemma~5 of \cite{LanguascoZaccagnini2010c}) yields
$\nu = 0.884472132 \dots$
Hence we can replace the value $-\log (0.91237810306)$ that we had in
\S7 of \cite{LanguascoZaccagnini2010c} with $-\log (0.884472132)$ in
the denominator of the definition of $s_0$ in the case where
$\lambda_1 / \lambda_2 \in \mathfrak{R}'$.


\providecommand{\bysame}{\leavevmode\hbox to3em{\hrulefill}\thinspace}
\providecommand{\MR}{\relax\ifhmode\unskip\space\fi MR }
\providecommand{\MRhref}[2]{%
  \href{http://www.ams.org/mathscinet-getitem?mr=#1}{#2}
}
\providecommand{\href}[2]{#2}

\smallskip
\author{Alessandro LANGUASCO\\
Universit\`a di Padova\\
Dipartimento di Matematica\\
Via Trieste 63\\
35121 Padova, Italy\\
E-mail: languasco@math.unipd.it}

\medskip
\author{Alessandro ZACCAGNINI \\
Universit\`a di Parma \\
Dipartimento di Matematica \\
Parco Area delle Scienze, 53/a \\
Campus Universitario \\
43124 Parma, Italy \\
E-mail: alessandro.zaccagnini@unipr.it}
\end{document}